\documentclass[12pt]{article}
\usepackage{amssymb}
\usepackage{amsmath}
\usepackage[latin1]{inputenc}
\usepackage[T1]{fontenc}
\usepackage{times}
\usepackage{dsfont}
\usepackage{natbib}
\usepackage{lineno}
\usepackage{mathrsfs}
\usepackage{hyperref}
\usepackage[english]{babel}
\input{bb.sty}
\date{}
\DeclareMathOperator{\Var}{Var}
\usepackage{graphicx}
\usepackage{xcolor}
\newtheorem{theorem}{Theorem}
\newtheorem{lemma}{Lemma}

\newtheorem{e-proposition}{Proposition}

\newtheorem{e-definition}[theoreme]{D\'efinition\rm}

\renewcommand{\theequation}{\arabic{equation}}

\numberwithin{equation}{section}

\begin{document}

\centerline{A Berry-Esseen theorem for sample quantiles under association }

\vskip 8mm
\noindent Lahcen DOUGE

\noindent FSTG, University Cadi Ayyad, B.P. 549 Marrakech, Morocco

\noindent lahcen.douge@uca.ac.ma

\vskip 3mm

\begin{abstract}
In this paper, the uniformly asymptotic normality for sample quantiles of associated random variables is investigated  under some conditions on the decay of the covariances. We obtain the rate of normal approximation of order $O(n^{-1/2}\log^2 n)$  if the covariances decrease exponentially to $0$. The best rate is shown as $O(n^{-1/3})$ under a polynomial decay of the covariances. 
\end{abstract}
\vskip 4mm

\section{Introduction} 
The concept of association for real-valued random variables was introduced by \cite{esary1967} and has found many applications in reliability theory and statistics. Let $(X_n)_{n\geq1}$ be a sequence of real-valued random variables. It is said to be associated if, for every finite subcollection $X_{i_1},\ldots,X_{i_n}$ and every pair of coordinatewise nondecreasing functions $g$, $h$ : $\mathbb{R}^{n}\rightarrow   \mathbb{R}$,
\begin{eqnarray*}
	\mathrm{Cov}\left( g(X_{i_1},\ldots,X_{i_n}), h(X_{i_1},\ldots,X_{i_n})\right) \geq 0,
\end{eqnarray*}
whenever the covariance is defined.

Assume that $(X_n)_{n\geq1}$ is a sequence of stationary associated random variables with a common marginal distribution function $F$ and finite second moment. For $0< p <1$, let
\begin{eqnarray*}
	x_p:= F^{-1}(p)= \inf\{x, F(x)\geq p\}
\end{eqnarray*}
denote the $p$th quantile of $F$. An estimator of $F^{-1}(p)$, $0< p <1$, is given by the sample $p$th quantile \begin{eqnarray*}
	x_{n,p}:= F_n^{-1}(p)= \inf\{x, F_n(x)\geq p\},
\end{eqnarray*}
where $F_n(x)= n^{-1}\sum_{i=1}^{n} I(X_{i}\leq x)$, $x\in\mathbb{R}$, denotes the empirical distribution function of $X_{1},\ldots, X_{n}$ and where $I(A)$ denotes the indicator function of a set $A$.

In this paper, we establish a Berry-Esseen theorem for sample quantiles of associated random variables under both exponential and polynomial decay of the covariances. The rate of normal approximation is of order $O(n^{-1/2}\log^2 n)$ under the exponential decay of the covariances. Under the polynomial decay of the covariances, the best rate obtained is of order $O(n^{-1/3})$.
\cite{lahiri2009} proved a Berry-Esseen theorem for sample quantiles of $\alpha$-strongly-mixing random variables under a polynomial mixing rate $\alpha(n)=O(n^{-\beta})$ with $\beta>12$. Their result has an optimal rate of order $O(n^{-1/2})$. \cite{yang2012berry} obtained the rate of order $O(n^{-1/6} \log n)$ under the condition on the mixing coefficient $\alpha(n)=O(n^{-\beta})$ with $\beta>39/11$. For more works on the Berry-Esseen bounds of sample quantiles, one can refer to \cite{yang2014note}, \cite{yang2012}, \cite{szewczak2017berry} and \cite{wang2019}.

The best known rate in the central limit theorem for associated random variables was obtained by  \cite{birkel1988} under an exponentially decaying of covariances. He obtained  a rate of order $O\big(n^{-1/2}\log n\big)$. Under a power decay of the covariance,  \cite{louhichi2002} obtained a rate of order $O\big(n^{-1/3}\big)$. For other papers about Berry-Esseeen bound for associated random variables, we can refer to  \cite{wood1983},  \cite{bulinskii1996} and  \cite{cai1999}.

Some applications of this results can be found in insurance where many risks are associated with heavy-tailed distribution and quantile based methods play an important role in evaluation of management risks. We refer to  \cite{denuit2006} for more details.

The organization of the paper is as follows. The basic assumptions and main results are listed in Section 2. In Section 3, some preliminary lemmas are given. The proofs of the main results are provided in Sections 4 and 5.

\section{Assumptions and main results}
We shall make use of the following conditions:
\begin{enumerate}
	\item[(C.1)] (i) $X_1$ has a bounded probability density function $f$.\\
	(ii) In a neighborhood of $x_p$, $F$ possesses a bounded second derivative $F^{''}$.
	\item[(C.2)] There exist constants $a_0\geq 0$ and $a>0$ such that for all $k\geq 1$,
	$$\mathrm{Cov}\left( X_1, X_{k+1}\right) \leq a_0\exp(-ak).$$
	\item[(C.3)] There exist constants $a_1\geq 0$ and $\beta >6$ such that for all $k\geq 1$,
	$$\mathrm{Cov}\left( X_1, X_{k+1}\right) \leq a_1 k^{-\beta}.$$
\end{enumerate}
Condition (C.1)(i) allows to derive the following covariance inequality, established by \cite{yu1993} for associated random variables, which will be needed in the proofs. 
\begin{eqnarray}\label{eq21}
\sup_{x,y\in \mathbb{R}}	\mathrm{Cov}\big( I(X_1\geq x), I(X_{k}\geq y)\big)  \leq A_0 \mathrm{Cov}(X_1,X_k)^{1/3}, \quad k\geq 1,
\end{eqnarray}
where the constant $A_0$ depends on the bound of $f$.
(C.1)(ii)  is the condition assumed by \cite{serfling2009} to show  the Berry-Esseen bound of the sample quantiles for identically independent random variables.\\
Denote 
\begin{eqnarray*}
	\sigma^{2}(x_p):=\mathrm{Var}\big( I(X_1\leq x_p)\big)  + 2\sum_{j=2}^{\infty}   \mathrm{Cov}\big( I(X_1\leq x_p), I(X_{j}\leq x_p)\big).       
\end{eqnarray*}

\begin{theorem}\label{thm1} Suppose that conditions (C.1) and (C.2) hold. If $f(x_p)>0$, then there exists a constant $C$ not depending on $n$ such that, for all large enough $n$,
	\begin{eqnarray}\label{thm1eq0}
		\sup_{t\in \mathbb{R}} \left| P\left(\frac{\sqrt{n}\big(F_n^{-1}(p)-F^{-1}(p)\big)}{a_p}\leq t\right)-\Phi(t)\right| \leq C\, n^{-1/2}\log^{2} n, 
	\end{eqnarray}
	where $a_p:=\sigma(x_p)/f(x_p)$ and $\Phi$ is the distribution function of a standard normal variable.
\end{theorem}
\begin{theorem}\label{thm2} Suppose that conditions (C.1) and (C.3) hold. If $f(x_p)>0$, then there exists a constant $B$ not depending on $n$ such that, for all large enough $n$,
	\begin{eqnarray*}\label{thm2eq0}
		\sup_{t\in \mathbb{R}} \left| P\left(\frac{\sqrt{n}\big(F_n^{-1}(p)-F^{-1}(p)\big)}{a_p}\leq t\right)-\Phi(t)\right| \leq B\,   \big(n^{-q/(8+2q)}+n^{-1/3}\big),
	\end{eqnarray*}
	where $q=2\left[ (\beta -3)/3\right]$. $\left[ x\right]$ stands for the integral part of $x$. 
\end{theorem}
\noindent The best rate of convergence in the central limit theorem under the arithmetic decay of the covariances is  $O\big(n^{-1/3}\big)$ and is obtained when $\beta\geq 15$.  This rate is of order $O\big(n^{-1/6}\big)$ if $\beta > 6$.\\\\
\noindent {\large \textbf{Simulation study}}\\
We perform a simulation study in order to investigate the rates of convergence in Theorem \ref{thm1} and Theorem \ref{thm2}. We consider the random sequence
\begin{eqnarray*}
	X_k= \sum_{j=1}^{m} a_j \varepsilon_{k+j}, \quad k,\, m\geq 1,
\end{eqnarray*}
where $ (a_j)_{j\in \mathbb{N}}$ is a sequence of real positive numbers and $ (\varepsilon_j)_{j\in \mathbb{N}}$ is a sequence of independent and identically distributed random variables with distribution  $\mathcal{N}(0,\sigma^2)$, $\sigma^2= 1/\sum_{j=1}^{m} a_j^2$. 
\begin{table}[h!]
	\begin{center}
	\caption{The uniform Berry-Esseen bounds}\label{tab:table1}
		\vspace{0.2cm}  
		\begin{tabular}{|cc||cccc||c}
			\hline\hline
		& 	&$p$&$n=100$  &$n=500$ &$n=1000$ \\ 
			\hline \hline
			
			&case 1 :&0.2  &$0.03636$  &$0.030211$ &$0.022639$ \\ 
			
		& $a_j=0.1^j$ &0.7 &$0.031639$ &$0.025639$&$0.018788$ \\ 
			\hline\hline
			&case 2 :&0.2  &$0.07336$  &$0.059064$ &$0.044885$ \\ 
			
			& $a_j=j^{-7}$ &0.7 &$0.062082$ &$0.053345$&$0.042656$ \\ 
			\hline\hline
		\end{tabular} 
	\end{center}
\end{table}
$(X_k)_{k\geq 1}$ is a stationary sequence of associated random variables \big(cf. ($\mathcal{P}_4$) of \cite{esary1967}\big) and $X_k\sim  \mathcal{N}(0,1)$, $k\geq 1$.
 Take $m=100$ and generate some random samples of  $(X_k)$  with the sample size $n =100$, $500$ and $1000$. We consider both the cases of the exponential and polynomial decay of the covariances and we compute the uniform Berry Esseen bounds in (\ref{thm1eq0}) for each value of $p=0.2$ or $0.7$. We first compute the term $A_n(t):=\sqrt{n}\big(F_n^{-1}(p)-F^{-1}(p)\big)/a_p\leq t$
and then evaluate the maximum value of  $\left|P (A_n(t)) -\Phi(t)\right|$ for  $t \in [-4, 4]$. The results are shown in  Table \ref{tab:table1}. We remark that,  for each value of $p$, the Berry Esseen bounds decrease as $n$ increases. These simulations show a good agreement with our main results.
 \section{Preliminaries}
\noindent Define, for each $i\geq 1$ and $t\in \mathbb{R}$,
\begin{eqnarray*}
	Y_i(t):= I\big(X_i\leq x_p+t\, a_p \,n^{-1/2}\big)- E I\big(X_i\leq x_p+t \,a_p \,n^{-1/2}\big)
\end{eqnarray*}
and 
\begin{eqnarray*}
	\sigma_{(n,t)}^2:=\Var\big(Y_1(t)\big)+ 2 \sum_{j=2}^{\infty}\mathrm{Cov}\big(Y_1(t), Y_j(t)\big).
\end{eqnarray*}
\begin{lemma}
	Suppose that conditions (C.1) and (C.2) hold. Then, for all large enough $n$ and any $t\in \mathbb{R}$ with $|t| < L_n:=b_0\log n$, $b_0>0$,
	\begin{eqnarray}\label{eq311}
	\left|\sigma_{(n,t)}^2-\sigma^2(x_p)\right| \leq C_1\, n^{-1/2} \log^2 n,
	\end{eqnarray}
	where $C_1$ is a constant not depending on $n$ and $t$.
\end{lemma}

\noindent \textbf{Proof.} 
For each $t\in \mathbb{R}$, the function  $G(x)=-I(x \leq x_p+t \,a_p \,n^{-1/2})$, $x\in \mathbb{R}$, is increasing. Then, by Property 4 in \cite{esary1967}, the sequence $\big(Y_i(t)\big)_{i\geq 1}$ is associated. \\
Let $t\in \mathbb{R}$ such that $|t| < L_n$, $n \geq 2$. First, by condition (C.2) and  (\ref{eq21}), we get $\sigma^2(x_p)<\infty$. Next, by condition (C.1) and Taylor's expansion, we have
\begin{eqnarray}\label{lem2eq1}
\lefteqn{\Big| 	\Var(Y_1(t)) -	\Var\big(I(X_1\leq x_p)\big) \Big| }\nonumber\\
&=&\Big| F \big(x_p+t \,a_p \,n^{-1/2}\big) -F (x_p)+F^2 (x_p) -F^2 \big(x_p+t \,a_p \,n^{-1/2}\big) \Big| \nonumber\\
&=&\Big| F \big(x_p+t \,a_p \,n^{-1/2}\big) -F (x_p)\Big|\Big(1+\Big| F \big(x_p+t \,a_p \,n^{-1/2}\big)+F (x_p)\Big|\Big)\nonumber\\
&\leq& 3 f(x_p) \,|t| \,a_p \,n^{-1/2}+o\big(|t| \,a_p \,n^{-1/2}\big) \nonumber\\
&=& O\big(n^{-1/2} \log n\big).
\end{eqnarray}
Similarly, for any $j\geq 2$,
\begin{eqnarray}\label{lem2eq11}
\lefteqn{\Bigg| E\Big[ I\big(X_1\leq x_p+t \,a_p \,n^{-1/2}\big) I \big(X_j\leq x_p+t \,a_p \,n^{-1/2}\big)\Big]}\nonumber\\
&& \qquad \quad -	E\Big[ I\big(X_1\leq x_p+t \,a_p \,n^{-1/2}\big) I \big(X_j\leq x_p\big)\Big] \Bigg| \nonumber\\
&\leq& E\Big|  I \big(X_j\leq x_p+t \,a_p \,n^{-1/2}\big)-   I \big(X_j\leq x_p\big) \Big| \qquad  \qquad    \qquad  \qquad  \qquad   \qquad\nonumber\\
&\leq & 2 f(x_p) \,|t| \,a_p \,n^{-1/2}+o\big(|t| \,a_p \,n^{-1/2}\big)\nonumber \\ 
&=& O\big(n^{-1/2} \log n\big)
\end{eqnarray}
and
\begin{eqnarray}\label{lem2eq12}
\lefteqn{\Bigg| E\Big[ I\big(X_1\leq x_p+t \,a_p \,n^{-1/2}\big)\Big]E\Big[ I \big(X_j\leq x_p+t \,a_p \,n^{-1/2}\big)\Big]}\nonumber \\
&& \quad \qquad - E\Big[ I\big(X_1\leq x_p+t \,a_p \,n^{-1/2}\big)\Big] E\Big[I \big(X_j\leq x_p\big)\Big] \Bigg|\nonumber \\
&=&O\big(n^{-1/2} \log n\big).\qquad \qquad \qquad\qquad\qquad  \qquad  \qquad   \qquad \qquad  \qquad  \qquad   \qquad
\end{eqnarray}
Combining (\ref{lem2eq11}) with (\ref{lem2eq12}), we obtain, for every $j\geq 2$,
\begin{eqnarray}\label{lem2eq2}
\lefteqn{	\Big| \mathrm{Cov}(Y_1(t), Y_j(t)) -	\mathrm{Cov}\big(I(X_1\leq x_p), I(X_j\leq x_p)\big) \Big|} \nonumber\\
&\leq&  \,\Big| \mathrm{Cov}(Y_1(t), Y_j(t)) -	\mathrm{Cov}\big(Y_1(t), I(X_j\leq x_p)\big)\Big| \nonumber\\
&& +  \Big|\mathrm{Cov}\big(Y_1(t), I(X_j\leq x_p)\big) -	\mathrm{Cov}\big(I(X_1\leq x_p), I(X_j\leq x_p)\big) \Big|\nonumber\\
&=& O\big(n^{-1/2} \log n\big).
\end{eqnarray}
Then, by putting (\ref{lem2eq1}) and (\ref{lem2eq2}) together and using condition (C.2) and (\ref{eq21}), we obtain, for $a_2=3/(2 a)$ and all large enough $n$, 
\begin{eqnarray*}
	\left| \sigma_{(n,t)}^2-\sigma^2(x_p)\right|&\leq& \Big|\Var(Y_1(t)) -	\Var\big(I(X_1\leq x_p)\big) \Big|\\
	&& + 2 \sum_{j=2}^{[a_2 \log n]} 	\Big| \mathrm{Cov}(Y_1(t), Y_j(t)) -	\mathrm{Cov}\big(I(X_1\leq x_p), I(X_j\leq x_p)\big) \Big|\\
	&& + 2\sum_{j=[a_2\log n]+1}^{\infty} 	\big| \mathrm{Cov}(Y_1(t), Y_j(t))\big| \\
	&& + 2\sum_{j=[a_2\log n]+1}^{\infty} 	\Big| \mathrm{Cov}\big(I(X_1\leq x_p), I(X_j\leq x_p)\big) \Big|\\
	&=& O\big(n^{-1/2} \log^2 n\big).
\end{eqnarray*}
\cqfd\\
Since the sequence $\big( I(X_i\leq x_p)\big) _{i\geq 1}$ is associated and $F$ is continuous at $x_p$, $0<\sigma^{2}(x_p)$. Then, it follows from (\ref{eq311}) that $\sigma_{(n,t)}>0$ for $n$ large enough. 
\begin{lemma}\label{lem1}
	Suppose that conditions (C.1) and (C.2) hold. If $f(x_p)>0$, then, for each  $t\in \mathbb{R}$ with $|t| < L_n$ and for all large enough $n$,
	\begin{eqnarray*}\label{eq31}
	\sup_{x\in \mathbb{R}}	\left| P\left(\frac{\sum_{i=1}^{n} Y_i(t)}{\sqrt{n}\sigma_{(n,t)}}\leq x\right)-\Phi(x)\right| \leq C_2\left[  n^{-1/2} \log n + \frac{1}{n\sigma_{(n,t)}^2}\right],
	\end{eqnarray*}
	where
	$C_2$ is a positive constant not depending on $n$ and $t$.
\end{lemma}
\noindent \textbf{Proof.} Let $\sigma_n^2:=\Var\big(\sum_{i=1}^{n} Y_i(t)\big)$. By using stationarity, condition (C.2) and (\ref{eq21}), we have, for $n$ large enough,
\begin{eqnarray}\label{eq32}
\nonumber| \sigma_n^2-n \sigma_{(n,t)}^2| &=& \left| 2\sum_{j=2}^{n} (n-j+1)\mathrm{Cov}\big( Y_1(t),Y_j(t)\big)  - 2n\sum_{j=2}^{\infty} \mathrm{Cov}\big (Y_1(t),Y_j(t) \big) \right| \\\nonumber
&=&2 n\sum_{j=n+1}^{\infty} \mathrm{Cov}\big(Y_1(t),Y_j(t)\big) + 2\sum_{j=2}^{n} (j-1)\mathrm{Cov}\big( Y_1(t),Y_j(t)\big)  \leq C_3,\\
\end{eqnarray}
where $C_3$ is some positive constant not depending on $n$ and $t$. For any $x\in \mathbb{R}$ and $n$ large enough
\begin{eqnarray}\label{lem1eq3b}
\lefteqn{ \Bigg| P\left(\frac{\sum_{i=1}^{n} Y_i(t)}{\sqrt{n}\sigma_{(n,t)}}\leq x\right)-\Phi(x)\Bigg| 
}\nonumber\\
&&\leq  \Bigg| P\left(\frac{\sum_{i=1}^{n} Y_i(t)}{\sigma_n}\leq \frac{\sqrt{n}\sigma_{(n,t)}}{\sigma_n}x\right)-\Phi \left(\frac{\sqrt{n}\sigma_{(n,t)}}{\sigma_n}x\right)\Bigg| +
\Bigg| \Phi\left(\frac{\sqrt{n}\sigma_{(n,t)}}{\sigma_n}x\right)-\Phi(x)\Bigg| \nonumber\\
&&:= \Psi_1+\Psi_2.
\end{eqnarray}
By (\ref{eq311}) and (\ref{eq32}), we deduce that $\lim_{n\rightarrow \infty}\sigma_n^2/(n\sigma_{(n,t)}^2) =1$ and $\inf_{n\geq n_0} \sigma^2_n/n>0$ for $n_0$ large enough. Then, by using (\ref{eq21}) and by applying Lemma \ref{lemA2} in the appendix, we get $$\Psi_1\leq C_4\, n^{-1/2} \log n,$$
where $C_4$ is a positive constant not depending on $n$ and $t$. Now, it is easy to see that
\begin{eqnarray}\label{eq33}
\sup_{x\in \mathbb{R}} \big| \Phi(px)-\Phi(x)\big|& \leq& (2\pi e)^{-1/2}\left[ (p-1)I(p\geq 1) + \left(\frac{1}{p}-1\right)I(0<p< 1) \right] \nonumber\\
& \leq& (2\pi e)^{-1/2} \, \left|p-\frac{1}{p}\right|, \quad p>0.
\end{eqnarray}
Thus, using (\ref{eq32}) and (\ref{eq33}), we obtain, for some positive constant $C_5$,
\begin{eqnarray*}\label{lem1eq3}
	\Psi_2
	\leq \frac{C_5}{n\sigma_{(n,t)}^2},
\end{eqnarray*}
which completes the proof of the lemma.\cqfd

\section{Proof of Theorem \ref{thm1}}
Denote
\begin{eqnarray*}
	G_{n}(t)=P\left( \frac{\sqrt{n}(x_{n,p}-x_{p})}{a_{p}}\leq t \right).
\end{eqnarray*}
We will show that
\begin{eqnarray}\label{thm1eq1}
\sup_{|t|>L_{n} }\big| G_{n}(t)-\Phi(t)\big| = O(n^{-1/2})
\end{eqnarray}
and
\begin{eqnarray}\label{thm1eq2}
\sup_{|t|\leq L_{n} }\big| G_{n}(t)-\Phi(t)\big|= O\big(n^{-1/2} \log^2 n\big).
\end{eqnarray}
First consider (\ref{thm1eq1}). Note that
\begin{eqnarray}\label{thm1eq21}
\sup_{|t|>L_{n} }\big| G_{n}(t)-\Phi(t)\big| &=&\max\Big\lbrace \sup_{t>L_{n} }\big| G_{n}(t)-\Phi(t)\big| , \sup_{t<-L_{n} }\big| G_{n}(t)-\Phi(t)\big| \Big\rbrace \nonumber\\
&\leq& \max\Big\lbrace 1-G_{n}(L_{n})+1-\Phi(L_{n}), G_{n}(-L_{n})+\Phi(-L_{n}) \Big\rbrace\nonumber\\
&\leq&  G_{n}(-L_{n})+1-G_{n}(L_{n})+1-\Phi(L_{n})\nonumber\\ 
&\leq&  P\big( | x_{n,p}-x_{p}| \geq a_{p}L_{n}n^{-1/2} \big)+1-\Phi(L_{n}). 
\end{eqnarray}
It is easily to see that
\begin{eqnarray}\label{thm1eqb22}
	1-\Phi(x)\leq \frac{(2\pi)^{-1/2}}{x} e^{-x^{2}/2},\quad x>0,
\end{eqnarray}
so that
\begin{eqnarray}\label{thm1eq22}
1-\Phi(L_{n})=O(n^{-1/2}).
\end{eqnarray}
Let $b_n:=(a_p-\varepsilon_0) n^{-1/2} L_n$, where $0<\varepsilon_0<a_p$. Note that
\begin{eqnarray*}
	P\big( | x_{n,p}-x_{p}| \geq a_{p}L_{n}n^{-1/2} \big)\leq P\big( \big| x_{n,p}-x_{p}\big| > b_n \big)
\end{eqnarray*}
and
\begin{eqnarray}\label{eq3}
P\big( | x_{n,p}-x_{p}| > b_n \big)=  P\big( x_{n,p}>x_p+b_n \big)+ P\big( x_{n,p}< x_p-b_n \big).
\end{eqnarray}
By Lemma \ref{lemA7}	, we obtain 
\begin{eqnarray}\label{eq4}
\nonumber P\big( x_{n,p}> x_p+b_n \big) &= & P\big( p> F_n(x_{p}+b_n) \big)= P\big(1- F_n(x_{p}+b_n)> 1-p \big)\\\nonumber
&=&P\Bigg(\sum_{i=1}^{n} I(X_i > x_p+b_n)> n(1-p) \Bigg)\\
&=&P\left(\sum_{i=1}^{n} \big[ I(X_i > x_p+b_n)- E I(X_i > x_p+b_n)\big]  > n \delta_1 \right),
\end{eqnarray}
where $\delta_1=F(x_{p}+b_n)-p$. Likewise,
\begin{eqnarray}\label{eq5}
P\big( x_{n,p}< x_p-b_n \big) \leq P\left(\sum_{i=1}^{n} \big[I(X_i \leq x_p-b_n)- E I(X_i \leq x_p-b_n)\big] \geq n \delta_2 \right),
\end{eqnarray}
where $\delta_2=p-F(x_{p}-b_n)$. Since $F$ is continuous at $x_p$, $F(x_p)=p$. By Taylor's expansion, one has
\begin{eqnarray*}
	F(x_p+b_n)-p&=&f(x_p)b_n+\frac{1}{2}\, F^{"}(\xi_1)b_n^2,\quad x_p<\xi_1<x_p+b_n,
\end{eqnarray*}
and
\begin{eqnarray*}
	p-F(x_p-b_n)&=&f(x_p)b_n -\frac{1}{2} \,F^{"}(\xi_2)b_n^2,\quad x_p-b_n<\xi_2<x_p.
\end{eqnarray*}
Let $M$ be a constant  such that, for $n$ large enough,
\begin{eqnarray*}
	\sup_{|z| \leq b_n}\big|F^{"}(x_p+z)\big| \leq M<\infty.
\end{eqnarray*}
Then
\begin{eqnarray}\label{thm1eq22b}
\delta :=  \min \left\lbrace \delta_1,\, \delta_2  \right\rbrace \geq b_n  \big[ f(x_p)- \frac{1}{2} M b_n\big].
\end{eqnarray}
By (\ref{eq3})-(\ref{thm1eq22b}) and condition (C.2), we obtain, by using (\ref{eq21}) and Lemma \ref{lemA1}, for all large enough $n$ and a suitable choice of $b_0$,
\begin{eqnarray*}\label{thm1eq23}
	P\big( | x_{n,p}-x_{p}|> b_n \big)=O(n^{-1/2}),
\end{eqnarray*}
which, when combined with (\ref{thm1eq21}) and (\ref{thm1eq22}), completes the proof of (\ref{thm1eq1}).\\

Now, we turn to (\ref{thm1eq2}). By Lemma \ref{lemA7}, we obtain, for all large enough $n$,
\begin{eqnarray*}
	G_{n}(t)=P\left( x_{n,p}\leq x_{p} + t a_{p} n^{-1/2} \right)&=&P\left( p\leq F_{n}\left( x_{p}+t a_{p} n^{-1/2}\right) \right) \\
	&=&	P\left( np\leq \sum_{i=1}^{n} I\left( X_{i}\leq x_{p} + t a_{p} n^{-1/2}\right) \right) \\
	&=&	P\left(  \frac{ \sum_{i=1}^{n} Y_{i}(t)}{\sqrt{n}\sigma_{(n,t)}}\geq - b_n(t)\right), 
\end{eqnarray*}
where
\begin{eqnarray*}
	b_n(t) = \frac{\sqrt{n}\left(F\left( x_{p}+t a_{p} n^{-1/2}\right)-p\right)}{\sigma_{(n,t)}}.
\end{eqnarray*}
Thus, it is easy to check that
\begin{eqnarray}\label{thm1eq6}
\Phi (t)-G_{n}(t) =	P\left( \frac{ \sum_{i=1}^{n} Y_{i}(t)}{\sqrt{n}\sigma_{(n,t)}}< - b_n(t)\right) -\Phi(-b_n(t))+\Phi(t)- \Phi(b_n(t)).
\end{eqnarray}
Next, by Taylor's expansion, one has
\begin{eqnarray}\label{thm1eq11}
b_n(t) &=& \frac{\sqrt{n}}{\sigma_{(n,t)}}  \left(F\left( x_{p}+t a_{p} n^{-1/2}\right)-F(x_p)\right)\nonumber\\
&=& \frac{\sqrt{n}}{\sigma_{(n,t)}}  \left(f( x_{p})t a_{p} n^{-1/2}  +\frac{1}{2}F^{''}(\xi_p)\left(t a_{p} n^{-1/2} \right)^2\right) \nonumber\\
&=& t\frac{\sigma(x_p)}{\sigma_{(n,t)}} +t^2 \frac{a_p^2 F^{''}(\xi_p)}{2\sigma_{(n,t)}} n^{-1/2},
\end{eqnarray}
where $ x_p< \xi_p< x_{p}+t a_{p} n^{-1/2}$. For each $t\in \mathbb{R}$ with $|t| < L_n$, it follows from (\ref{eq311}) and (\ref{thm1eq11}) that $|b_n(t)| <\infty$ for $n$ large enough. Then, by Lemma \ref{lem1}, one has 
\begin{eqnarray*}
	\left| P\left( \frac{ \sum_{i=1}^{n} Y_{i}(t)}{\sqrt{n}\sigma_{(n,t)}}< - b_n(t)\right) -\Phi(-b_n(t))\right|
	& \leq & C_1\left[  n^{-1/2} \log n + \frac{1}{n\sigma^2_{(n,t)}}\right].
\end{eqnarray*}
By (\ref{eq311}), we deduce that, for all large enough $n$,
\begin{eqnarray}\label{thm1eq7}
\sup_{|t| \leq L_n}\left| P\left( \frac{ \sum_{i=1}^{n} Y_{i}(t)}{\sqrt{n}\sigma_{(n,t)}}< - b_n(t)\right) -\Phi(-b_n(t))\right| & \leq & M_1\, n^{-1/2} \log n,
\end{eqnarray}
where $M_1$ is a positive constant. Now, for $n$ large enough,
\begin{eqnarray}\label{thm1eq9}
\sup_{|t| \leq L_n}\left| \Phi(t) -\Phi(b_n(t))\right| \leq \sup_{|t| \leq L_n}\left| \Phi(t) -\Phi\left( \frac{\sigma(x_p)}{\sigma_{(n,t)}}t\right)\right|+ \sup_{|t| \leq L_n}\left| \Phi\left( \frac{\sigma(x_p)}{\sigma_{(n,t)}}t\right) -\Phi(b_n(t))\right|.\nonumber\\
\end{eqnarray}
By (\ref{eq33}) and (\ref{eq311}), we get that for $n$ large enough and for some constant $M_2$,
\begin{eqnarray}\label{thm1eq10}
\nonumber\sup_{|t| \leq L_n}\left| \Phi(t) -\Phi\left( \frac{\sigma(x_p)}{\sigma_{(n,t)}}t\right)\right|&\leq& \frac{M_2}{\sigma^2(x_p)} \sup_{|t| \leq L_n}\left| \sigma_{(n,t)}^2-\sigma^2(x_p)\right| \\
&=&O\left( n^{-1/2} \log^2 n\right).
\end{eqnarray}
Now, observe that
\begin{eqnarray}\label{thm1eq12}
\sup_{x\in \mathbb{R}}\left| \Phi\left(x+y\right) -\Phi(x)\right|\leq | y|  (2\pi)^{-1/2},\quad \forall y \in \mathbb{R}.
\end{eqnarray}
Hence, by condition C1(ii) with (\ref{thm1eq11}) and (\ref{thm1eq12}), we obtain for all large enough $n$,
\begin{eqnarray}\label{thm1eq13}
\sup_{|t| \leq L_n}\left| \Phi\left( \frac{\sigma(x_p)}{\sigma_{(n,t)}}t\right) -\Phi(b_n(t))\right| = O\left( n^{-1/2}\log^2 n\right).
\end{eqnarray}
Finally, it follows from (\ref{thm1eq6}), (\ref{thm1eq7})-(\ref{thm1eq10}) and (\ref{thm1eq13}) that
\begin{eqnarray*}\label{thm1eq14}
	\sup_{|t| \leq L_n}\left| G_n(t)-\Phi(t)\right| =O\left( n^{-1/2}\log^2 n\right).
\end{eqnarray*}
This completes the proof of Theorem \ref{thm1}. \cqfd
\section{Proof of Theorem \ref{thm2}}
Let $K_n:=n^{b}$, $b=1/(4+q)$. We will show that
\begin{eqnarray}\label{thm2eq01}
\sup_{|t|>K_{n} }\big| G_{n}(t)-\Phi(t)\big| = O\left( n^{-q/(4+q)}\right)
\end{eqnarray}
and
\begin{eqnarray}\label{thm2eq02}
\sup_{|t|\leq K_{n} }\big| G_{n}(t)-\Phi(t)\big| =O\big(n^{-q/(8+2q)}+n^{-1/3}\big).
\end{eqnarray}
Let $c_n:=(a_p-\varepsilon_1) n^{b-1/2} $, where $0<\varepsilon_1<a_p$. From (\ref{eq3})-(\ref{eq5}), we can assert that
\begin{eqnarray}\label{thm2eq1}
P\big( | x_{n,p}-x_{p}| \geq a_{p}K_{n}n^{-1/2}\big)&\leq& P\Big( \big| x_{n,p}-x_{p}\big| > c_n \Big)\nonumber\\
 &\leq& P\left(\sum_{i=1}^{n} V_i\geq n \alpha_1 \right)+ P\left(\sum_{i=1}^{n} W_i \geq n \alpha_2 \right),
\end{eqnarray}
where $V_i=I(X_i > x_p+c_n)- E I(X_i > x_p+c_n)$, $W_i=I(X_i \leq x_p-c_n)- E I(X_i \leq x_p-c_n)$, $i\geq 1$, $\alpha_1=F(x_{p}+c_n)-p$ and $\alpha_2=p-F(x_{p}-c_n)$.\\

As in (\ref{thm1eq22b}), we have for $n$ large enough,
\begin{eqnarray*}
	\alpha :=  \min \left\lbrace \alpha_1,\, \alpha_2  \right\rbrace \geq c_n  \big[ f(x_p)- \frac{1}{2} M c_n\big]>M_3\, n^{b-1/2},
\end{eqnarray*}
where $M_3>0$ is some constant. $(V_i)_{i\geq 1}$ and  $(W_i)_{i\geq 1}$ are two bounded associated sequences. Then, by  (\ref{eq21}) and Lemma \ref{lemA6}, we get for any $r\geq 1$,
\begin{eqnarray*}\
	\sup \left| \mathrm{Cov}(V_{t_1}\ldots V_{t_m} , V_{t_{m}+r}\ldots V_{t_q})\right|&\leq& \sup
	\sum_{i=t_1}^{t_m}\sum_{j=t_m+r}^{t_q}   \mathrm{Cov}\big(I(X_i \geq x_p+c_n), I(X_j\geq x_p+c_n)\big)\\
	&\leq& M_4 \,\sup\sum_{i=t_1}^{t_m}\sum_{j=t_m+r}^{t_q}   \big( \mathrm{Cov}(X_i,X_j)\big)^{1/3}\\
	&=& O(r^{-q/2}),
\end{eqnarray*}
where $M_4$ is some positive constant and  $1\leq t_1\leq \ldots \leq t_m \leq t_{m}+r\leq \ldots \leq t_q$. Then, one can apply Lemma \ref{lemA4} and Markov inequality to show that
\begin{eqnarray}\label{thm2eq3}
P\left(\sum_{i=1}^{n} V_i \geq n \alpha_1 \right) &\leq & n^{-bq} n^{-q/2} E\left[ \sum_{i=1}^{n} V_i\right]^{q}\nonumber\\
&=& O\left( n^{-bq}\right).
\end{eqnarray}
Likewise 
\begin{eqnarray}\label{thm2eq4}
P\left(\sum_{i=1}^{n} W_i \geq n \alpha_2 \right) = O\left( n^{-bq}\right).
\end{eqnarray}
Thus, from (\ref{thm2eq1})-(\ref{thm2eq4}), we get
\begin{eqnarray}\label{thm2eq5}
P\big( | x_{n,p}-x_{p}| \geq a_{p}K_{n}n^{-1/2}\big)=O\left( n^{-bq}\right).
\end{eqnarray}
The proof of (\ref{thm2eq01}) is deduced from (\ref{thm1eq21}), (\ref{thm1eqb22}) and (\ref{thm2eq5}).\\

Now, let us establish (\ref{thm2eq02}). By the same arguments as in (\ref{lem2eq1}) and (\ref{lem2eq2}), we obtain, for any $t\in \mathbb{R}$ with $|t| < K_n$,
\begin{eqnarray*}
	\left| \sigma_{(n,t)}^2-\sigma^2(x_p)\right|&\leq& \Big|\Var(Y_1(t)) -	\Var\big(I(X_1\leq x_p)\big) \Big|\nonumber\\
	&& +2 \sum_{j=2}^{[ n^{b}]} 	\Big| \mathrm{Cov}(Y_1(t), Y_j(t)) -	\mathrm{Cov}\big(I(X_1\leq x_p), I(X_j\leq x_p)\big) \Big|\nonumber\\
	&&  + 2\sum_{j=[n^{b}]+1}^{\infty} 	\big| \mathrm{Cov}(Y_1(t), Y_j(t))\big| + 2\sum_{j=[n^{b}]+1}^{\infty} 	\Big| \mathrm{Cov}\big(I(X_1\leq x_p), I(X_j\leq x_p)\big) \Big|\nonumber\\
	&=& O \left[ n^{b-1/2}+ n^{2b-1/2}+ n^{-bq/2}\right], \nonumber\\
\end{eqnarray*}
which implies, for any $t\in \mathbb{R}$ with $|t| < K_n$,
\begin{eqnarray}\label{thm2eq6}
\left| \sigma_{(n,t)}^2-\sigma^2(x_p)\right|&=& O\big(n^{-q/(8+2q)}\big).
\end{eqnarray}
By condition (C.3), similar to the proof of (\ref{eq32}), we get , for $n$ large enough, 
\begin{eqnarray}\label{thm2eq6b}
| \sigma_n^2-n \sigma_{(n,t)}^2| \leq M_5,
\end{eqnarray}
where $M_5>0$ is some constant not depending on $t$. On the other hand, we can check easily that, for stationary associated sequences, condition (i) in Lemma \ref{lemA3} is satisfied as soon as $\inf_{m\geq 1} \sigma_m^2/m>0$ which is verified by using (\ref{thm2eq6}), (\ref{thm2eq6b}) and the fact that $\mathrm{Var}(Y_1(t))>0$ for $n$ large enough.
Therefore,  by using Lemma \ref{lemA3} with $\delta=1$, and the assertions (\ref{lem1eq3b}) , (\ref{eq33}), we obtain, for all large enough $n$,
\begin{eqnarray}\label{thm2eq7}
\sup_{|t| \leq K_n}\left| P\left( \frac{ \sum_{i=1}^{n} Y_{i}(t)}{\sqrt{n}\sigma_{(n,t)}}< - b_n(t)\right) -\Phi(-b_n(t))\right|
& \leq & O\left[  n^{-1/3} + \frac{1}{n\sigma^2_{(n,t)}}\right]\nonumber\\
& \leq & O \big( n^{-1/3}\big) .
\end{eqnarray}
Next, by arguing as in the proof of (\ref{eq33}) and by using (\ref{thm2eq6}), we have 
\begin{eqnarray*}\label{thm2eq8}
	\sup_{|t| \leq K_n}\left| \Phi(t) -\Phi\left( \frac{\sigma(x_p)}{\sigma_{(n,t)}}t\right)\right|\leq \frac{M_6}{\sigma^2(x_p)} \sup_{|t| \leq K_n}\left| \sigma_{(n,t)}^2-\sigma^2(x_p)\right| =O\big(n^{-q/(8+2q)}\big),
\end{eqnarray*}
where $M_6>0$. Hence, by (\ref{thm1eq11}) and (\ref{thm1eq12}), we get 
\begin{eqnarray*}\label{thm2eq9}
	\sup_{|t| \leq K_n}\left| \Phi\left( \frac{\sigma(x_p)}{\sigma_{(n,t)}}t\right) -\Phi(b_n(t))\right| = O\left( n^{-q/(8+2q)}\right).
\end{eqnarray*}
Thus
\begin{eqnarray}\label{thm2eq10}
\sup_{|t| \leq K_n}\left| \Phi(t) -\Phi(b_n(t))\right| &\leq& \sup_{|t| \leq K_n}\left| \Phi(t) -\Phi\left( \frac{\sigma(x_p)}{\sigma_{(n,t)}}t\right)\right|+ \sup_{|t| \leq K_n}\left| \Phi\left( \frac{\sigma(x_p)}{\sigma_{(n,t)}}t\right) -\Phi(b_n(t))\right|\nonumber\\
&=&  O\big(n^{-q/(8+2q)}\big).
\end{eqnarray}
Finally, using (\ref{thm1eq6}), (\ref{thm2eq7}) and (\ref{thm2eq10}), we obtain
\begin{eqnarray*}\label{thm2eq11}
	\sup_{|t| \leq K_n}\left| G_n(t)-\Phi(t)\right| =O\big(n^{-q/(8+2q)}+n^{-1/3}\big).
\end{eqnarray*}
This completes the proof of Theorem \ref{thm2}. \cqfd

\newpage

\begin{appendix}
	\section{Appendix}
	\setcounter{lemma}{0}
	\renewcommand{\thelemma}{\Alph{section}\arabic{lemma}}
	\begin{lemma} [ \cite{douge2007}, Theorem 1] \label{lemA1}
		Let $(X_i)_{i\in\mathbb{N}}$ be  a stationary sequence of associated random variables such that
		\begin{enumerate}
			\item[i)]$\vert X_{i}\vert\le A_1$, $\forall~i\ge 0$, where $A_1$ is
			some constant.
			\item[ii)]$\mathrm{Cov}(X_1, X_{k+1})\leq \theta_{0}
			\exp(-\theta k)$, for any $k\geq 0,~ \theta> 0$ and
			$\theta_{0}>0$.
		\end{enumerate}
		Then, for all $n\ge 2$ and all
		$\varepsilon>\frac{6 A_1}{\sqrt{n}}$, one has
		\begin{eqnarray*}\label{A1}
			P\left(\frac{1}{n}\left|\sum_{i=1}^n (X_i-\mathbb{E}X_i)\right|\ge
			\varepsilon\right)\leq
			8 A_{2}\exp\left(-\frac{\theta\wedge1}{12 A_1}\sqrt{n}\,\varepsilon\right),
		\end{eqnarray*}
		where
		$A_{2}=\exp\left(\dfrac{\theta_{0}}{4 A_1^{2}(1-e^{-\theta})}\right).$
	\end{lemma}
	
	\begin{lemma}[ \cite{birkel1988}, Theorem 2.1] \label{lemA2}
		Let $(X_i)_{j\geq 1}$ be an associated process with $EX_j=0$, $j\geq 1$, satisfying
		\begin{enumerate}
			\item[i)]    $u(n):= \sup_{k\in \mathbb{N}}\sum_{j: | j-k| \geq n} \mathrm{Cov}(X_j, X_k)\leq c_0 \exp (-c_1 n)$ for some $c_0>0$ and $c_1>0$,
			\item[ii)]   $\exists n_0 \in \mathbb{N}^{*}$ such that $\inf_{n\geq n_0} \sigma_n^2/n>0$, where $\sigma_n^2:=\Var\big(\sum_{i=1}^{n} X_i\big)$\\
			\noindent and
			\item[iii)]  $\sup_{j\in \mathbb{N}} E|X_j|^{3+c_2}<\infty$ for some $c_2>0$.
		\end{enumerate}
		Then there exists a constant $A_3$ not depending on $n$ such that for any $n\geq n_0$
		\begin{eqnarray*}\label{A2}
			\sup_{x\in \mathbb{R}} \left| P\left(\frac{\sum_{i=1}^{n} X_i}{\sigma_n}\leq x\right)-\Phi(x)\right| \leq A_3\, n^{-1/2} \log n.
		\end{eqnarray*}
	\end{lemma}
	
	\begin{lemma}[ \cite{louhichi2002}, Theorem 1] \label{lemA3}
		Let $(X_i)_{i\in\mathbb{N}}$ be  a stationary sequence of centered and bounded associated random variables with second finite moment such that
		\begin{enumerate}
			\item[i)]
			\begin{eqnarray*}\label{A3}
				\exists n_0 \in \mathbb{N}^{*}\quad \text{such that}\quad \inf_{n\geq n_0,\, 0\leq k<n} \frac{\sigma_n^2-\sigma_k^2}{n-k}>0.
			\end{eqnarray*}
			\item[ii)]
			\begin{eqnarray*}\
				\sum_{i=1}^{\infty} i^{\delta}\mathrm{Cov}(X_1, X_{1+i})\textcolor{red}{<}\infty, \quad 0<\delta \leq 1.
			\end{eqnarray*}
		\end{enumerate}
		Then,  for any $n\geq n_0$,
		\begin{eqnarray*}\label{A2}
			\sup_{x\in \mathbb{R}} \left| P\left(\frac{\sum_{i=1}^{n} X_i}{\sigma_n}\leq x\right)-\Phi(x)\right| = O( n^{-\delta/3}).
		\end{eqnarray*}
	\end{lemma}	
	
	\begin{lemma}[ \cite{doukhan1999}, Theorem 1]\label{lemA4}
		Let $(X_i)_{i\in\mathbb{N}}$ be a sequence of random variables fulfilling for some fixed $q\in \mathbb{N}$, $q\geq 2$
		\begin{eqnarray*}\
			\sup \left| \mathrm{Cov}(X_{t_1}\ldots X_{t_m} , X_{t_{m+1}}\ldots X_{t_q})\right| =O(r^{-q/2})\quad \text{as}\, r \rightarrow \infty,
		\end{eqnarray*}
		where the supremum is taken over all $\left\lbrace t_1,\ldots, t_q\right\rbrace $ such that $1\leq t_1\leq \ldots \leq t_q$ and $m$, $r$ satisfy $t_{m+1}-t_{m}=r$. Then there exists a constant $A_4$ not depending on $n$ for which
		\begin{eqnarray*}
			\left| E S_n^q\right|  \leq A_4 n^{q/2}.
		\end{eqnarray*}
	\end{lemma}	
	
\begin{lemma}[ \cite{birkel1988}, Lemma 3.1]\label{lemA6}	
		Let $A$ and $B$ be finite sets and let $X_j$, $j\in A\cup B$, be associated random variables. Then for all  real-valued partially differentiable functions $g$, $h$ with bounded partial derivatives, there holds
		 	\begin{eqnarray*}
			\left| \mathrm{Cov}\left( g\left( (X_i)_{i\in A}\right) , h\left( (X_j)_{j\in B}\right) \right)  \right|  \leq \sum_{i\in A}\sum_{j\in B} \left\| \partial g/\partial t_i\right\| _{\infty}\left\| \partial h/\partial t_j\right\| _{\infty}  \mathrm{Cov}(X_i,Y_j).
		\end{eqnarray*}
	\end{lemma}	
	
	\begin{lemma}[ \cite{serfling2009}, Lemma 1.1.4]\label{lemA7}	
		Let $F$ is the distribution function. The function $F^{-1}(t)$, $0<t<1$, is non decreasing and left continuous, and satisfies
		$$F(x)\geq t \quad \text{if and only if} \quad x\geq F^{-1}(t).$$
		\end{lemma}	
\end{appendix}

\bibliographystyle{natbib}
\bibliography{Berry-Esseen}

\end{document}